\documentstyle[11pt]{article}
\newcommand{\bbz}{\mbox{\boldmath $Z$}}

\newcommand{\bbc}{\mbox{\boldmath $C$}}

\newcommand{\bbf}{{F}}
\newcommand{\qed}{\hfill$\Box$}

\newenvironment{note}[1]{\par\addvspace{\medskipamount}\noindent
                         {\bf {#1}}\sl
                       }{\par\addvspace{\medskipamount}\rm}

\newcommand{\suml}{\sum\limits}
\newcommand{\de}{\partial}
\newcommand{\half}{{1\over 2}}
\font\gotic=eufm10
\newcommand{\gots}{{\hbox{\gotic\char'123}}}

\begin{document}

\title{Hecke Algebra Actions on the Coinvariant Algebra}
\bibliographystyle{acm}
\author{Ron M. Adin%
\thanks{Department of Mathematics and Computer Science, Bar-Ilan University,
Ramat Gan 52900, Israel. Email: {\tt radin@math.biu.ac.il} }$^{\ \S}$
\and
Alexander Postnikov%
\thanks{Department of Applied Mathematics, Massachusetts
Institute of Technology, MA 02139, USA. Email: {\tt apost@math.mit.edu} }
\and Yuval Roichman%
\thanks{Department of Mathematics and Computer Science, Bar-Ilan University,
Ramat Gan 52900, Israel. Email: {\tt yuvalr@math.biu.ac.il} } 
\thanks{Research supported in part by  the Israel Science Foundation
and by internal research grants from Bar-Ilan University.}}
\date{(submitted: August 1, 1999; revised: April 13, 2000)}

\newcommand{\TT}{{T_1, \ldots, T_k}}
\newcommand{\lamlam}{{\lambda^1, \ldots, \lambda^k}}

\maketitle
\begin {abstract}

Two actions of the Hecke algebra of type A on the corresponding
polynomial ring are studied.  Both are deformations of 
the natural action of the symmetric group on polynomials, 
and keep symmetric functions invariant.
We give an explicit description of these actions, and deduce a
combinatorial formula for the resulting graded characters
on the coinvariant algebra.

\end{abstract}

\section{\bf Introduction}

\noindent{\bf 1.1.} 
The symmetric group $S_n$ acts on the polynomial
ring $P_n=\bbf[x_1,\dots,x_n]$ (where $\bbf$ is a field of characteristic 
zero) by permuting variables.
Let $I_n$ be the ideal of $P_n$ generated by the symmetric ($S_n$-invariant)
polynomials without a constant term.
The {\it coinvariant algebra} of type $A$ is the quotient $P_n/I_n$.
Schubert polynomials, constructed in the seminal papers [BGG] and [De], 
form a distinguished basis for the coinvariant algebra. These
polynomials correspond to Schubert cells in the corresponding flag variety.

\smallskip

\noindent{\bf 1.2.}
In this paper we present two deformations of this action.
For these deformations we can take $\bbf=\bbc(q)$, the field of rational 
functions in an indeterminate $q$. Most of the results actually hold when 
$\bbf$ is replaced by the ring $\bbz[q]$ of polynomials in $q$ with integer 
coefficients.

\smallskip

Let $T_1,\dots,T_{n-1}$ be the standard generators of the Hecke algebra
${\cal H}_n(q)$ of type $A$; for definitions see Section 2.1 below.

The first action
$\rho_1:{\cal H}_n(q) \to Hom_{\bbf}(P_n,P_n)$
is defined using $q$-commu\-ta\-tors:
$$
\rho_1(T_i):=\de_iX_i-qX_i\de_i\qquad (1\le i<n), \leqno (1.1)
$$
where 
$$
\de_i:={1\over x_i-x_{i+1}}(1-s_i)
$$ 
is the divided difference operator (see Section 2.2),
and $X_i$ denotes multiplication by $x_i$.
This action belongs to a family introduced in [LS]
(see Subsection 7.1 below). For a geometric interpretation see [DKLLST].
In [DKLLST, Section 1] such families of operators are attributed to 
Hirzebruch~[Hr]. 

\smallskip

The second action is naturally defined on monomials by the formula
$$
\rho_2(T_i) (x_i^\alpha x_{i+1}^\beta m):=\cases
{q x_i^\beta x_{i+1}^\alpha m, & \hbox{ if $\alpha> \beta$;} \cr
(1-q) x_i^\alpha x_{i+1}^\beta m+  x_i^\beta x_{i+1}^\alpha m, 
  & \hbox{ if $\alpha< \beta$;} \cr
x_i^\alpha x_{i+1}^\beta m, & \hbox{ if $\alpha= \beta$.}} \leqno (1.2)
$$
Here $m$ is a monomial involving neither $x_i$ nor $x_{i+1}$. 

For a closely related action (defined in the context of quantum groups)
 see [Ji].

\medskip

\begin{note}
\noindent{\bf Claim:} 
The ideal $I_n$ is invariant under both actions.
The resulting graded characters on the coinvariant algebra
have a common combinatorial formula. 
\end{note}

This shows, in particular, that $\rho_1$ and $\rho_2$ lead to equivalent 
representations of ${\cal H}_n(q)$ on the coinvariant algebra. 
For $q=1$ they both reduce to the natural $S_n$ action.

\medskip

\noindent{\bf 1.3.} 
Since the ideal $I_n$ is invariant under both $\rho_1$ and $\rho_2$, the
coinvariant algebra $P_n/I_n$ carries appropriate actions $\tilde{\rho}_1$
and $\tilde{\rho}_2$.  Let $\chi_1^k$ and $\chi_2^k$ be the characters of 
these representations on the $k$-th homogeneous component of $P_n/I_n$.
We shall give an explicit formula for these characters, using the following 
combinatorial function.

For any permutation $w\in S_n$, define
$$
m_q(w):=\cases
{(-q)^m,& \hbox{ if there exists a unique $0\le m< n$ so that}\cr
& \hbox{ $w(1)>\dots>w(m+1)<w(m+2)<\dots<w(n)$;}\cr
0, & \hbox { otherwise. }\cr} \leqno (1.3)
$$
Let $\mu:=(\mu_1,\dots,\mu_t)$ be a partition of $n$,
and let $S_{\mu} := S_{\mu_1} \times \ldots \times S_{\mu_t}$ 
be the corresponding Young subgroup of $S_n$.
For any permutation $w\in S_n$ write
$w=r\cdot (w_1\times\cdots\times w_t)$, where
$w_i\in S_{\mu_i}$ $(1\le i\le t)$ and 
$r$ is a representative of minimal length for the left coset 
$w S_\mu$ in $S_n$.  Define
$$
\hbox{weight}_q^{\mu}(w):=\prod\limits_{i=1}^t m_q(w_i). \leqno (1.4)
$$

\smallskip

\begin{note}
\noindent{\bf Theorem.} For all $k\ge 0$ and $\mu\vdash n$,
$$
\chi_1^k(T_\mu)=\chi_2^k(T_\mu)=
\suml_{\{w\in S_n:\ell(w)=k\}} \hbox{weight}_q^\mu(w) 
$$
where $T_\mu:=T_1 T_2 \cdots T_{\mu_1-1} T_{\mu_1+1} \cdots\cdots
T_{\mu_1+\ldots+\mu_t-1}$ is the subproduct of $T_1 T_2 \cdots T_{n-1}$
omitting $T_{\mu_1+\dots+\mu_i}$ for all $1\le i<t$.
\end{note}

\medskip

The proof relies on an explicit description of the action with respect to 
the Schubert basis of the coinvariant algebra.
See Theorem 4.1 and Theorem 6.5 below.

\smallskip

\noindent{\bf Remark.} 
This character formula is a natural $q$-analogue of a weight formula for $S_n$
presented in [Ro2].
A formally similar result appears also in Kazhdan-Lusztig theory. 
Kazhdan-Lusztig characters may be represented as sums of 
exactly the same weights, but over different summation sets 
[Ro1 Corollary 4, Ra2]. 

\medskip

\noindent{\bf 1.4.} 
The rest of this paper is organized as follows.
Preliminaries and necessary background are given in Section 2.
In Section 3 we introduce $q$-commutators and study their 
representation matrices.
The character formula for $q$-commutators is proved in Section 4.
Natural randomized operators are introduced in Section 5.
In Section 6 we show that the representations induced by the two different
actions are equivalent. 
Section 7 concludes the paper with remarks regarding 
related families of operators, connections with Kazhdan-Lusztig theory,
and open problems. 

\bigskip

\section{\bf Preliminaries}

\medskip

\subsection{\bf The Hecke Algebra of Type $A$}

\medskip

The symmetric group $S_n$ is generated by $n-1$ involutions
$s_1,s_2,\dots,s_{n-1}$ satisfying the Moore-Coxeter relations
$$
s_is_{i+1}s_i=s_{i+1}s_is_{i+1} 
 \hbox{ }\hbox{ }\hbox{ }\hbox{ } (1\le i< n-1)\leqno (2.1)
$$
and
$$
s_is_j=s_js_i  \hbox{ }\hbox{ }\hbox{ }\hbox{ if }  |i-j|>1. \leqno (2.2)
$$
These involutions are known as the Coxeter generators of $S_n$.

All reduced expressions of a permutation $w\in S_n$ 
with respect to these generators have the same
length, denoted by $\ell(w)$.

The Hecke algebra ${\cal H}_n(q)$ of type $A$ is the algebra over 
$\bbf := \bbc(q)$
generated by $n-1$ generators  $T_1,\dots,T_{n-1}$, 
satisfying the Moore-Coxeter relations (2.1) and (2.2) as well as 
the following ``deformed involution'' relation:
$$
T_i^2= (1-q)T_i+q \hbox{ }\hbox{ }\hbox{ }\hbox{ } (1\le i<n).
\leqno (2.3)
$$

It should be noted that the last relation is slightly 
non-standard; this is done in order to get more elegant
$q$-analogues.
 In order to shift to the standard version,
one should replace each $T_i$ by $-T_i$.

Let $w$ be a permutation in $S_n$ and let
$s_{i_1}\cdots s_{i_{\ell(w)}}$ be a reduced expression for $w$.
It follows from the above relations that
$T_w:= T_{i_1}\cdots T_{i_{\ell(w)}}$ is independent
of the choice of reduced expression;
the set $\{T_w| w\in S_n\}$ forms a linear basis for
${\cal H}_n(q)$.

Let $\mu=(\mu_1,\dots,\mu_t)$ be a partition of $n$.
Define $T_\mu\in {\cal H}_n(q)$ to be the product 
$$
T_\mu:=T_1T_2\cdots T_{\mu_1-1}T_{\mu_1+1}T_{\mu_1+2}
\cdots T_{\mu_1+\mu_2-1}T_{\mu_1+\mu_2+1}\cdots\cdots T_{\mu_1+\dots+\mu_t-1}.
$$
This is the subproduct of the product
$T_1T_2\cdots T_{n-1}$ of all generators (in the usual order), 
obtained by omitting $T_{\mu_1+\dots+\mu_i}$ for all $1\le i< t$.
These elements play an important role in the character theory of 
${\cal H}_n(q)$.
For $q=1$, the elements $T_\mu$ are representatives
of all conjugacy classes in $S_n$. 
It follows that, for $q=1$, a character is determined by its values 
at these elements. This is also the case for arbitrary $q$,
as the following theorem shows.

\smallskip

\begin{note}
\noindent{\bf Theorem 2.1.} [Ra1, Theorem 5.1]
For each $w\in S_n$ there exists a linear combination 
$$
C_w= \suml_{\mu} a_{w\mu} T_\mu\in {\cal H}_n(q),
$$
with $a_{w\mu}\in \bbz[q]$, such that 
$$
\chi(T_w)=\chi(C_w)
$$
for all characters $\chi$ of the Hecke algebra ${\cal H}_n(q)$.
\end{note}

\medskip

Let $\mu=(\mu_1,\dots,\mu_t)$ be a partition of $n$.
Each permutation $w\in S_n$ has an associated  weight, 
$\hbox{weight}_q^\mu(w)$, as defined in (1.3)--(1.4).
The irreducible characters of ${\cal H}_n(q)$ are indexed by the
partitions of $n$. These characters may be represented
as weighted sums over Knuth equivalence classes.

\smallskip

\begin{note}
\noindent{\bf Theorem 2.2.} [Ro1, Corollary 4] 
Let ${\cal C}$ be a Knuth equivalence class of shape $\lambda$. Then
$$
\chi^\lambda(T_\mu)=\suml_{w\in \cal C} \hbox{weight}_q^\mu(w),
$$
where $\chi^\lambda$ is the irreducible character of ${\cal H}_n(q)$ corresponding 
to the shape $\lambda$.
\end{note}

\bigskip

\subsection{\bf  Schubert Polynomials and the Coinvariant Algebra}

\medskip

\subsubsection{\bf Basic Actions on the Polynomial Ring}

Let $x_1,x_2,\dots,x_n$ be independent variables, and let $P_n$
be the polynomial ring $\bbf[x_1,x_2,\dots,x_n]$.
The symmetric group $S_n$ acts on $P_n$ by permuting the 
variables $x_i$.
Let $\Lambda_n=\Lambda[x_1,x_2,\dots,x_n]$ be the subring of symmetric functions
(i.e., polynomials which are invariant under the action of $S_n$).
Denote  by $\Lambda_n(i)$ the ring of all polynomials
which are invariant under the action of $s_i$ for a fixed $i$, $1\le i<n$.
Clearly,  $f\in \Lambda_n(i)$ if and only if $f$ is
symmetric in the variables $x_i$ and $x_{i+1}$. We call 
the polynomials in $\Lambda_n(i)$, {\it $i$-symmetric} polynomials.

\medskip

For $1\le i<n$ define a divided difference operator  $\partial_i:P_n\rightarrow P_n$ by
$$
\de_i:=(x_i-x_{i+1})^{-1}(1-s_i).
$$
If $f\in P_n$ is a homogeneous polynomial of degree $d$, 
which is not $i$-symmetric, then $\de_i (f)$
is homogeneous of degree $d-1$.
For $i$-symmetric polynomials $\partial_i (f)=0$.

The operators $\de_i$ satisfy the nil-Coxeter relations [Ma, (2.1)]:
$$
 \de_i^2=0 \hbox{ }\hbox{ }\hbox{ }\hbox{ } (1\le i<n),\leqno (2.4)
$$
$$
\de_i\de_{i+1}\de_i=\de_{i+1} \de_i \de_{i+1} 
\hbox{ }\hbox{ }\hbox{ }\hbox{ } (1\le i<n-1),\leqno (2.5)
$$
$$
\de_i \de_j= \de_j \de_i 
   \hbox{ }\hbox{ }\hbox{ if $|i-j|>1$}. \leqno (2.6)
$$

Let $X_i$ be the operator on $P_n$ corresponding to
multiplication by $x_i$.
Clearly, $X_i$ increases degree by $1$. 

The algebra generated by the operators $\de_i, 1\le i<n$, and
 $X_i, 1\le i\le n$ was studied in [De] and [BGG]. 
The generators
satisfy the following commutation relations:

$$
\de_i X_j=X_j\de_i
   \hbox{ }\hbox{ }\hbox{ if $|i-j|>1$}, \leqno (2.7)
$$
$$
\de_i X_i = 1+ X_{i+1}\de_i 
\hbox{ }\hbox{ }\hbox{ }\hbox{ } (1\le i<n),\leqno (2.8)
$$
$$
X_i \de_i = 1+ \de_i X_{i+1} 
\hbox{ }\hbox{ }\hbox{ }\hbox{ } (1\le i<n),\leqno (2.9)
$$

\medskip

\subsubsection{\bf Schubert Polynomials}

For any sequence  $a=(a_1,\dots,a_k)$ of positive integers less than $n$,
 define $\de_a:=\de_{a_1}\cdots\de_{a_k}$.
It follows from the relations (2.5)--(2.6) that if $a$, $b$ are two reduced 
expressions for the same permutation $w\in S_n$ then $\de_a=\de_b$.
We can therefore use the notation $\de_w$ for $w\in S_n$, and in particular
$\de_{s_i}:=\de_i$ for $1\le i <n$.

The relation $\de_i^2=0$ implies that for any $w\in S_n$
and any $1\le i< n$

$$
\de_i\de_w=\cases
{\de_{s_i w}, & \hbox{ if $\ell(s_i w)>\ell(w)$;} \cr
0, & \hbox{ if $\ell(s_i w)<\ell(w)$.} \cr} \leqno (2.10)
$$

\smallskip

For each $w\in S_n$ we define the $Schubert$ $polynomial$ $\gots_w$
by
$$
\gots_w:=\de_{w^{-1}w_0}(x_1^{n-1}x_2^{n-2}\cdots x_{n-1})
$$
where $w_0$ is the longest element in $S_n$.

By definition, $\gots_w$ is a homogeneous polynomial of degree $\ell(w)$.

It follows from (2.10) that 
$$
\de_i (\gots_w)=\cases
{\gots_{ws_i}, & \hbox{ if $\ell(ws_i)<\ell(w)$;} \cr
0, & \hbox{ if $\ell(ws_i)>\ell(w)$.} \cr} \leqno (2.11)
$$
Denote $\gots_{s_i}$ by $\gots_i$.
For any $1\le i< n$,
$$
\gots_i=x_1+\cdots+x_i. \leqno (2.12)
$$
See [Ma, (4.4)]. The following is an important variant of Monk's formula.
 
\smallskip

\begin{note}
\noindent{\bf Monk's Formula.} [Ma, (4.11)] 
Let $1\le i< n$ and $w\in S_n$.  Then
$$
\gots_i\gots_w=\sum\limits_t \gots_{wt}
$$
where the sum extends over all transpositions $t=t_{jk}$ 
interchanging $j$ and $k$,
with $1\le j\le i<k\le n$ and $\ell(wt)=\ell(w)+1$.
\end{note}

\medskip

The description of the action of the operator $X_i$ on Schubert polynomials 
follows from Monk's formula and (2.12):
$$
X_i(\gots_w) = (\gots_i - \gots_{i-1})\gots_w =
\sum_{j\le i< k} \gots_{wt_{jk}} - \sum_{j< i\le k} \gots_{wt_{jk}}
\leqno (2.13)
$$
$$
= \sum_{j=i<k} \gots_{wt_{jk}} - \sum_{j<i=k} \gots_{wt_{jk}},
$$
where all summations are over the transpositions $t=t_{jk}$ satisfying 
$\ell(wt) = \ell(w) + 1$, with $j$ and $k$ in the indicated ranges.

\smallskip

\subsubsection{\bf The Coinvariant Algebra}

Recall that $\Lambda_n=\Lambda[x_1,\dots,x_n]$ is the subring of $P_n$
consisting of symmetric functions, and
let $I_n$ be the ideal of $P_n$ generated by symmetric functions without a 
constant term.
The quotient $P_n/I_n$ is called the $coinvariant$ $algebra$ of $S_n$.
$S_n$ acts naturally on this algebra. The resulting representation 
is isomorphic to the regular representation of the symmetric group.
See, e.g., [Hu, \S 3.6] and [Hi, \S II.3].

Let $R^k$ ($0\le k \le {n\choose 2}$) be the $k$-th homogeneous component of 
the coinvariant algebra : $P_n/I_n=\oplus_{k=0}^{n\choose 2} R^k$. 
Each $R^k$ is an $\bbf[S_n]$-module; 
let $\chi^k$ be the corresponding character. The set $\{\gots_w|w\in S_n\}$
of Schubert polynomials forms a basis for $P_n/I_n$, and 
the set $\{\gots_w|\ell(w)=k\}$ forms a basis for $R^k$.

\smallskip

The action of the simple reflections on  Schubert polynomials
is described by the following proposition, which is a reformulation
of [BGG, Theorem 3.14 (iii)].

\smallskip

\begin{note}
\noindent{\bf Proposition 2.3.} For any simple reflection $s_i$ 
and any $w\in S_n$,
$$
s_i(\gots_w)=\cases
{\gots_w, & \hbox { if $\ell(ws_i)>\ell(w)$;}\cr
- \gots_w 
+ \suml_{k<i} \gots_{w(k,i+1,i)}
- \suml_{k<i} \gots_{w(k,i,i+1)} +\cr
+ \suml_{k>i+1} \gots_{w(k,i,i+1)}
- \suml_{k>i+1} \gots_{w(k,i+1,i)},
 &\hbox{ if $\ell(ws_i)<\ell(w)$,}\cr}
$$
where  $(k,i,i+1),(k,i+1,i)$ are cycles of length 3,
and the sums extend over those values of $k$ (in the prescribed
ranges) for which $w(k,i,i+1)$ (respectively, $w(k,i+1,i)$ ) has the same
length as $w$.
\end{note}

\smallskip

Note that the signs in this proposition may depend on notational conventions.

Let $\mu$ be a partition of $n$, and let $\chi^k$ be the $S_n$-character 
on $R^k$ as above.
The following character formula is analogous to Theorem 2.2.

\smallskip

\begin{note}
\noindent{\bf Theorem 2.4.} [Ro2, Theorem 2] 
With the notations of Theorem 2.2,
$$
\chi^k(w_\mu)=\suml_{\ell(w)=k} \hbox{weight}_1^\mu(w)
$$
where $ \hbox{weight}_1^\mu(w)$ is the weight (1.4) with $q=1$,
and $w_\mu$ is any permutation of cycle-type $\mu$.
\end{note}

\smallskip

The goal of this paper is to define a Hecke algebra action
on the polynomial ring $P_n$ which produces a $q$-analogue of Theorem 2.4.

\bigskip

\section{\bf $q$-Commutators}

\medskip

For $1\le i< n$ define the $q$-commutator $[\de_i,X_i]_q$ as follows :
$$
[\de_i,X_i]_q:=\de_i X_i-q X_i \de_i.
$$
It should be noted that for $q=1$, $[\de_i,X_i]_1=s_i$.
Let $A_i:=[\de_i,X_i]_q$.

\smallskip

\begin{note}
\noindent{\bf Claim 3.1.} The operators $A_i$, $1\le i< n$, satisfy 
the Hecke algebra relations (2.1)--(2.3).
\end{note} 

\noindent{\bf Proof.} 
Combine the 
 nil-Coxeter relations (2.4)--(2.6) for the operators $\de_i$
with the commutation relations (2.7)--(2.9) for the operators
$\de_i$  and $X_j$.

\qed

\smallskip

It follows that the mapping $T_i\mapsto A_i$ ($1\le i<n$) may be extended to a representation $\rho_1$
of ${\cal H}_n(q)$ on $P_n$:
$$
\rho_1(T_i):=A_i=[\de_i,X_i]_q.
$$

\smallskip

\noindent{\bf Remark.} The polynomial action of the Coxeter generators 
of $S_n$ is {\it multiplicative}, i.e.,
for any generator $s_i$ and any two polynomials $f,g\in P_n$: 
$$
s_i (fg)= s_i(f) s_i(g). \leqno (3.1)
$$
Thus each $s_i$ acts on $P_n$ as an {\it algebra} automorphism. 
 It follows that if $f$ is $i$-symmetric (see Section 2.2.1)
 then  $s_i (fg)= f s_i(g)$.
In contrast to that, the operators $A_i$ are not multiplicative.
Actually, (2.3) implies that the eigenvalues of any linear action of 
a Hecke algebra generator $T_i$ are $1$ and $-q$,
and taking $f$ to be a $(-q)$-eigenvector of $A_i$,
one would get (if $A_i$ were multiplicative) that $f^2$ is a $q^2$-eigenvector,
which is impossible for generic $q$.

\medskip

\begin{note}
\noindent{\bf Claim 3.2.} For any $1\le i<n$, 
any $i$-symmetric polynomial $f\in \Lambda_n(i)$, 
and any polynomial $g\in P_n$: 
$$
A_i (f)=f \leqno (3.2)
$$
and
$$
A_i (fg)=f A_i (g) \leqno (3.3)
$$
\end{note}

\noindent{\bf Proof.} (3.2) is the special case $g=1$ of (3.3). The latter follows from the fact that for arbitrary polynomials $f,g\in P_n$
$$
\de_i(fg)=\de_i(f)g+s_i(f)\de_i(g).
$$
Therefore, if $f\in\Lambda_n(i)$ then $\de_i(fg)=f\de_i(g)$, so that 
$$
A_i(fg)=\de_i(x_ifg)-qx_i\de_i(fg)=f[\de_i(x_ig)-qx_i\de_i(g)]=fA_i(g).
$$
\qed

\smallskip

It follows that the ideal $I_n$ of $P_n$ is invariant under all the operators $A_i$, giving rise to a representation $\tilde \rho_1$ of ${\cal H}_n(q)$ on the quotient $P_n/I_n$, namely:
on the coinvariant algebra.
Let $\chi^k_1$ be the character of this representation on
the $k$-th homogeneous component $R^k$ of $P_n/I_n$ ($0\le k\le {n \choose 2}$).

Recall from Section 2.2.3 that the set of Schubert polynomials
$\{\gots_w|\ell(w)=k\}$ forms a basis for $R^k$.

The representation $\tilde \rho_1$ yields a $q$-analogue of Proposition 2.3. 

\smallskip

\begin{note}
\noindent{\bf Theorem 3.3.} 
For any $1\le i< n$ and $w\in S_n$,
$$
\tilde\rho_1(T_i) (\gots_w)=\cases{
\gots_w, & \hbox { if $\ell(ws_i)>\ell(w)$;}\cr
- q\gots_w 
+ q\suml_{k<i} \gots_{w(k,i+1,i)}
-  \suml_{k<i} \gots_{w(k,i,i+1)} +\cr
+  \suml_{k>i+1} \gots_{w(k,i,i+1)}
- q\suml_{k>i+1} \gots_{w(k,i+1,i)},
&\hbox { if $\ell(ws_i)<\ell(w)$,}\cr
}
$$
where  $(k,i,i+1),(k,i+1,i)$ are cycles of length 3,
and the sums extend over those values of $k$ (in the prescribed
ranges) for which $w(k,i,i+1)$ (respectively, $w(k,i+1,i)$ ) has the same
length as $w$.
\end{note}

\smallskip

\noindent{\bf Proof.} By the commutation relation (2.8),
$$
A_i=1+(X_{i+1}-qX_{i})\de_i.
$$
Applying  (2.11)
and (2.13) completes the proof.
\qed

\bigskip

\section{\bf Characters of $q$-Commutators}

In this section we prove the following $q$-analogue of Theorem 2.4.

\smallskip

\begin{note}

\noindent{\bf Theorem 4.1.} For any partition $\mu\vdash n$ and $k\ge 0$,
$$
\chi_1^k(T_\mu)=\suml_{\ell(w)=k} \hbox{weight}_q^\mu(w)
$$
where $\hbox{weight}_q^\mu(w)$ is defined as in $(1.4)$, and
the subproduct $T_\mu$ is defined as in Section 2.1.
\end{note}

\medskip

First recall that, by Theorem 3.3, for any $1\le i<n$ and $w\in S_n$
$$
\tilde A_i(\gots_w)=\cases{
\gots_w,  &\hbox { if $\ell(ws_i)>\ell(w)$;}\cr
-q\gots_w + \suml_{\ell(zs_i)>\ell(z)=\ell(w)} a_{w,z}(q) \gots_{z},
&\hbox { if $\ell(ws_i)<\ell(w)$,}\cr
} \leqno (4.1)
$$
where $\tilde A_i:=\tilde\rho_1(T_i)$, $a_{w,z}(q)\in \bbz[q]$ and 
the summation is over all $z\in S_n$ with $\ell(zs_i)>\ell(z)=\ell(w)$.

\medskip

Denote by $\langle \cdot ,\cdot \rangle$ the inner product on $P_n/I_n$ 
defined by $\langle \gots_v,\gots_w\rangle:=\delta_{vw}$,
where $\delta_{vw}$ is the Kronecker delta. In order to prove Theorem 4.1 we need the following lemma.

\begin{note}

\noindent{\bf Lemma 4.2.} 
Let $w\in S_n$ be a permutation satisfying
$\ell(ws_i)<\ell(w)$. Then, for any $\pi\in S_n$:
$$
\langle \tilde A_i \tilde A_\pi (\gots_w), \gots_w\rangle= 
-q\langle \tilde A_\pi (\gots_w), \gots_w\rangle.
$$

\end{note}

\noindent{\bf Proof of Lemma 4.2.}
 It follows from $(4.1)$ that if $\ell(ws_i)<\ell(w)$ and $v\in S_n$ then
$$
\langle \tilde A_i (\gots_v),\gots_w\rangle=\cases{
-q,& if $v=w$;\cr
0,& if $v\not=w$.\cr} \leqno (4.2)
$$
Substituting (4.2) into 
\begin{eqnarray*}
\langle \tilde A_i \tilde A_\pi (\gots_w), \gots_w\rangle &=& 
\langle \tilde A_i (\suml_v \langle \tilde A_\pi(\gots_w),\gots_v\rangle
 \gots_v), \gots_w\rangle \cr
&=& \suml_v \langle \tilde A_\pi(\gots_w), \gots_v\rangle 
    \langle \tilde A_i(\gots_v), \gots_w\rangle
\end{eqnarray*}
we obtain the desired conclusion.

\qed

\medskip

\noindent{\bf Proof of Theorem 4.1.}

In order to prove Theorem 4.1, it suffices to prove that for any partition $\mu=(\mu_1,\dots,\mu_t)$ of $n$
$$
\langle \tilde A_\mu (\gots_w),\gots_w\rangle=\hbox{weight}_q^\mu(w),
$$
where $\tilde A_\mu=\tilde\rho_1(T_\mu)$ is the subproduct of
$\tilde A_1\tilde A_2\cdots\tilde A_{n-1}$ obtained by omitting $\tilde A_{\mu_1+\dots+\mu_i}$ for all $1\le i< t$.

Assume now that there is an index $i$ such that: 
$\tilde A_i$ and $\tilde A_{i+1}$ are factors of $\tilde A_\mu$, 
$\ell(ws_i)>\ell(w)$, and $\ell(ws_{i+1})<\ell(w)$. 
Then, by Lemma 4.2 :
$$
\langle \tilde A_{i+1} \tilde A_\mu (\gots_w), \gots_w\rangle=
-q\langle \tilde A_\mu (\gots_w), \gots_w\rangle.\leqno (4.3)
$$
On the other hand, by the Hecke algebra relations: 
$\tilde A_{i+1}\tilde A_\mu=\tilde A_\mu \tilde A_i$.
Hence
$$
\langle \tilde A_{i+1} \tilde A_\mu (\gots_w), \gots_w\rangle=
\langle  \tilde A_\mu \tilde A_i (\gots_w), \gots_w\rangle
=\langle  \tilde A_\mu (\gots_w), \gots_w\rangle. \leqno(4.4)
$$
The last equality follows from (4.1).

Comparing (4.3) and (4.4) we obtain 
$$
-q\langle  \tilde A_\mu(\gots_w), \gots_w\rangle= 
\langle  \tilde A_\mu(\gots_w), \gots_w\rangle.
$$
We conclude that, if  
there is an index $i$ so that $\tilde A_i$ and $\tilde A_{i+1}$ are factors
of $\tilde A_\mu$, $\ell(ws_i)>\ell(w)$, and $\ell(ws_{i+1})<\ell(w)$, 
then (since $q$ is indeterminate)
$$
\langle \tilde A_\mu (\gots_w), \gots_w\rangle=0.
$$
Note that in this case $i$, $i+1$ and $i+2$ belong to the same ``block'' in 
the partition $\mu$, and $w(i) < w(i+1) > w(i+2)$.  Thus indeed 
$$\hbox{weight}_q^{\mu}(w) = 0.$$

It remains to check the case in which 
there is no index $i$ so that both $\tilde A_i$ and $\tilde A_{i+1}$ appear as 
factors in the product $\tilde A_\mu$, 
with $\ell(ws_i)>\ell(w)$ and $\ell(ws_{i+1})<\ell(w)$. 

In this case, the relation $\tilde A_i\tilde A_j=\tilde A_j\tilde A_i$ for $|i-j|>1$ gives 
$$
\tilde A_\mu=\tilde A_{i_1} \cdots \tilde A_{i_m} \tilde A_{i_{m+1}}\cdots \tilde A_{i_{\mu_1+\dots+\mu_t-t}},
$$
where $\ell(ws_{i_j})<\ell(w)$ for $j\le m$, and $\ell(ws_{i_j})>\ell(w)$ for $j>m$.  
Applying (4.1) and Lemma 4.2 iteratively implies
$$
\langle \tilde A_\mu (\gots_w), \gots_w\rangle=(-q)^m=\hbox{weight}_q^\mu(w),
$$
where 
$m=\#\{i\,|\,\ell(ws_i)<\ell(w) 
             \hbox{ and $\tilde A_i$ is a factor of $\tilde A_\mu$}\}$.
\qed

\bigskip

\section{\bf Randomized Operators}

\medskip

In this section we define a natural
``randomized" action of the Coxeter generators
on the polynomial ring $P_n$, and show that this action satisfies
the Hecke algebra relations.
This action will be defined initially  on monomials, and
then  extended  by linearity to all polynomials in $P_n$.

Let $e_{\alpha,\beta,m}:=x_i^\alpha x_{i+1}^\beta m$,
where $m\in P_n$ is a monomial involving neither $x_i$
nor $x_{i+1}$, and $\alpha,\beta$ are nonnegative integers.
Note that the linear subspace 
$V_{\alpha,\beta,m}:={\rm span}\{e_{\alpha,\beta,m},e_{\beta,\alpha,m}\}$
is invariant under the action of $s_i$.
In this space $s_i$ acts as a transposition of the two basis elements
(if $\alpha \ne \beta$).

A natural randomization of $e_{\alpha,\beta,m}$ is
$(1-q)e_{\alpha,\beta,m}+q e_{\beta,\alpha,m}$, 
where the parameter $q$ may be interpreted as transition probability
$0\le q\le 1$.
Motivated by well-known asymmetric physical processes (simulated annealing etc.), we define
$$
R_i^* (e_{\alpha,\beta,m}):=\cases{
e_{\beta,\alpha,m}, & \hbox{ if $\alpha\ge \beta$;} \cr
(1-q)e_{\alpha,\beta,m}+q e_{\beta,\alpha,m}, & \hbox{ if $\alpha< \beta$,} \cr
} \leqno (5.1)
$$
and extend this randomized action to the whole polynomial ring $P_n$ by 
linearity. See also [Ji].

\medskip

\begin{note}
\noindent{\bf Claim 5.1.} The operators $R^*_i$, $1\le i< n$, satisfy 
the Hecke algebra relations (2.1)--(2.3).
\end{note}

\smallskip

\noindent {\bf Proof.} Easily verified by an explicit calculation of the action on the monomials $e_{\alpha,\beta,m}$.
\qed

\medskip

The operators $R^*_i$ lead, therefore, to a representation of 
${\cal H}_n(q)$ on $P_n$. Unfortunately, the symmetric functions are not invariant under this action.
Consider, therefore, the operators whose representation matrices with respect to the basis of monomials
are the transposes of those representing $R^*_i$; i.e., define
$$
R_i (e_{\alpha,\beta,m}):=\cases{
q e_{\beta,\alpha,m}, & \hbox{ if $\alpha> \beta$;} \cr
(1-q) e_{\alpha,\beta,m}+ e_{\beta,\alpha,m}, & \hbox{ if $\alpha< \beta$;} \cr
e_{\alpha,\beta,m}, & \hbox{ if $\alpha= \beta$.} \cr
} \leqno (5.2)
$$
Of course, the operators $R_i$, $1\le i<n$, also satisfy the Hecke relations 
(2.1)--(2.3).
It follows that the mapping $T_i\mapsto R_i$ ($1\le i<n$)
may be extended to a representation $\rho_2$ of ${\cal H}_n(q)$ on $P_n$:
$$
\rho_2(T_i):=R_i .
$$

The following claim is analogous to Claim 3.2.

\begin{note}
\noindent
{\bf Claim 5.2.} 
For any $1\le i<n$, 
any $i$-symmetric polynomial $f\in \Lambda_n(i)$, 
and any polynomial $g\in P_n$: 
$$
R_i (f)=f \leqno (5.3)
$$
and
$$
R_i (fg)=f R_i (g) \leqno (5.4)
$$
\end{note}

\noindent{\bf Proof.} Direct calculation. 
\qed

\medskip

It follows from the first part of the claim that 
symmetric functions are pointwise invariant under $\rho_2({\cal H}_n(q))$.
By the second part, the ideal $I_n$ is also invariant under 
$\rho_2({\cal H}_n(q))$.
Thus, $\rho_2$ gives rise to a representation $\tilde \rho_2$ of ${\cal H}_n(q)$ on the coinvariant algebra $P_n/I_n$.

\medskip

The action of $R_i$ on monomials is transparent. Section 6 is devoted to
a better understanding of the action on the coinvariant algebra.

\bigskip

\section{\bf Properties of the Randomized Action}

The following sequence of assertions concerns the connections between the operators $A_i$ and $R_i$.

\smallskip

\begin{note}
\noindent{\bf Claim 6.1.} The operators $A_i$ and $R_i$ have the same invariant vectors: 
$$
\ker(A_i-1) =\ker(R_i-1)=\Lambda_n(i),
$$
where $\Lambda_n(i)$ is the set (actually, subalgebra)
of all polynomials invariant under $s_i$.
\end{note}

\noindent{\bf Proof.}
By the definition of $A_i$ and the commutation relations~(2.8), 
$$
\ker(A_i-1)= \ker[(X_{i+1}-qX_i)\de_i]=\ker\de_i =\Lambda_n(i).
$$
As for $R_i-1$, let $V_{\alpha,\beta,m}:={\rm span}\{e_{\alpha,\beta,m},
e_{\beta,\alpha,m}\}$ as in the beginning of Section~5.
Note that
$$
P_n= \bigoplus_{\{(\alpha,\beta,m)|\alpha\ge \beta\}} V_{\alpha,\beta,m}
$$
is a decomposition of $P_n$ into a direct sum of $R_i$-invariant subspaces.

By (5.2), in $V_{\alpha,\beta,m}$ :
$$
(R_i-1)(e_{\alpha,\beta,m})=\cases{
-e_{\alpha,\beta,m}+qe_{\beta,\alpha,m},& if $\alpha>\beta$;\cr
-qe_{\alpha,\beta,m}+e_{\beta,\alpha,m},& if $\alpha<\beta$;\cr
0, & if $\alpha=\beta$.\cr}
$$
Thus
$$
\ker(R_i-1)\cap V_{\alpha,\beta,m}=\cases{
{\rm span}\{e_{\alpha,\beta,m}+e_{\beta,\alpha,m}\},& if $\alpha\ne\beta$;\cr
{\rm span}\{e_{\alpha,\beta,m}\},& if $\alpha=\beta$,\cr}
$$
implying
$$
\ker(R_i-1)=\bigoplus_{\{(\alpha,\beta,m)|\alpha\ge\beta\}}
{\rm span}\{e_{\alpha,\beta,m} + e_{\beta,\alpha,m}\}=\Lambda_n(i).
$$
\qed

\smallskip

\begin{note}
\noindent{\bf Claim 6.2.}
\begin{itemize}
\item[(a)] For any positive integers $i<n$, $j\le n$
 and nonnegative integer $m$,
$$
(A_i-R_i) (x_j^m)=(1-q)\de_i (x_j^{m+1}). \leqno (6.1)
$$
\item[(b)] For any polynomial $f \in \bbz[x_1,\ldots,x_n]$,
the polynomial $(A_i-R_i)f$ is $i$-symmetric and divisible by $1-q$ :
$$
(A_i-R_i)f \in \Lambda_n(i) \cap (1-q)\cdot\bbz[x_1,\ldots,x_n].
$$
\end{itemize}
\end{note}

\smallskip

\noindent{\bf Proof.}\\
(a) If $j\not \in \{i,i+1\}$ or $m=0$ then both sides of (6.1) equal zero.
If $j=i$ and $m\ge 1$ then
\begin{eqnarray*}
(A_i-R_i)(x_i^m) &=& (1+(x_{i+1}-qx_i)\de_i-R_i)(x_i^m) \cr
&=& x_i^m+(x_{i+1}-qx_i)\suml_{t=1}^m x_i^{m-t}x_{i+1}^{t-1}-qx_{i+1}^m \cr
&=& (1-q)\suml_{t=0}^m x_i^{m-t} x_{i+1}^t =(1-q)\de_i(x_i^{m+1}).
\end{eqnarray*}
If $j=i+1$ and $m\ge 1$ then, by Claim 6.1
\begin{eqnarray*}
(A_i-R_i)(x_{i+1}^m) &=& (A_i-R_i)(x_i^m+x_{i+1}^m-x_i^m)
= (A_i-R_i)(-x_i^m) \cr
&=& -(1-q)\de_i (x_i^{m+1})=(1-q)\de_i (x_{i+1}^{m+1}).
\end{eqnarray*}
(b) It suffices to prove this claim for monomials 
$x_1^{k_1}x_2^{k_2}\cdots x_n^{k_n}$.
Any such monomial has the form $gx_j^m$, where $g \in \Lambda_n(i)$, 
$j\in\{i,i+1\}$, and $m$ is a nonnegative
integer. It follows from (3.3), (5.4) and (6.1) that
$$
(A_i-R_i)(gx_j^m)=g(A_i-R_i)(x_j^m)=(1-q)g\de_i (x_j^{m+1}),\leqno (6.2)
$$
as claimed.
\qed

\smallskip

\begin{note}
\noindent{\bf Lemma 6.3.} $\Lambda_n(i)$ is spanned, as a $\Lambda_n$-module,
by the Schubert polynomials $\gots_w$ with $\ell(ws_i)>\ell(w)$. 
The same holds when the ground field $\bbf$ is replaced by $\bbz$.
\end{note}

\smallskip

\noindent{\bf Proof.} 
%
%
%
First of all, if $\ell(ws_i) > \ell(w)$ then, by Proposition~2.3, 
$s_i(\gots_w) = \gots_w$ and therefore $\gots_w \in \Lambda_n(i)$.

By the same proposition, for any $w\in S_n$
$$
(1+s_i)(\gots_w) \in {\rm span}\{\gots_z \,|\, \ell(zs_i)>\ell(z)=\ell(w)\}
\qquad(\mbox{\rm in\ } P_n/I_n)
$$
and therefore, for any $f\in P_n/I_n$ :
$$
(1+s_i)(f) \in {\rm span}\{\gots_z \,|\, \ell(zs_i)>\ell(z)\}
\qquad(\mbox{\rm in\ } P_n/I_n).
$$
If $f\in \Lambda_n(i)/I_n$ then $(1+s_i)(f) = 2f$, so that $\Lambda_n(i)/I_n$ 
is spanned, as a vector space, by the above Schubert polynomials.
Since $I_n$ is the ideal of $P_n$ generated by the homogeneous elements in
$\Lambda_n$ of positive degree, a standard argument yields the claimed result 
for $\Lambda_n(i)$ as a $\Lambda_n$-module.

The result for $\bbz$ instead of $\bbf$ follows from the fact that 
Schubert polynomials also form a ``basis'' for polynomials with integer 
coefficients.
\qed

\smallskip

The following proposition provides a description
of the action of $\tilde\rho_2({\cal H}_n(q))$ on Schubert polynomials.

\smallskip

\begin{note}
\noindent{\bf Proposition 6.4.} 
For each $1\le i<n$ and $w\in S_n$,
$$
\tilde \rho_2(T_i)(\gots_w)=\cases{
\gots_w, & \hbox { if $\ell(ws_i)>\ell(w)$;}\cr
-q\gots_w + \suml_{\ell(zs_i)>\ell(z)} 
[(1-q)b_{w,z}+c_{w,z}]\gots_z,
&\hbox { if $\ell(ws_i)<\ell(w)$,}\cr
}
$$
where $b_{w,z}\in \bbz$, $c_{w,z}\in \{-1,0,1\}$, 
and the sum extends over all permutations $z\in S_n$ with
$\ell(zs_i)>\ell(z)=\ell(w)$.
\end{note}

\medskip

\noindent{\bf Proof.}
Since $\gots_w\in \Lambda_n(i)$ for $w\in S_n$ with $\ell(ws_i>\ell(w)$, 
Claim 6.1 implies that
$\rho_2(T_i)(\gots_w)=\gots_w$ for these $w$.

Homogeneous components of $P_n$ are invariant
under the action of each $R_i, 1\le i<n$. 
It follows that the homogeneous components of the coinvariant algebra
are invariant under $\tilde\rho_2({\cal H}_n(q))$, so that each
$\tilde \rho_2(T_i)(\gots_w)$ is spanned by Schubert polynomials
of degree $\ell(w)$. Combining this fact with Claim 6.2(b) and Lemma 6.3 
shows that for any $1\le i<n$ and $w\in S_n$
$$
(\tilde\rho_2(T_i)-\tilde\rho_1(T_i))(\gots_w)=
(1-q)\suml_{\ell(zs_i)>\ell(z)=\ell(w)}d_{w,z} \gots_z,\leqno (6.3)
$$
where $d_{w,z}\in \bbz$,
and the sum extends over all permutations $z\in S_n$ with
$\ell(zs_i)>\ell(z)=\ell(w)$.

Combining (6.3) with (4.1) gives, for any $w\in S_n$ with $\ell(ws_i)<\ell(w)$ 
\begin{eqnarray*}
\tilde\rho_2(T_i)(\gots_w) 
&=&\tilde\rho_1(T_i)(\gots_w)+(\tilde\rho_2(T_i)-\tilde\rho_1(T_i))(\gots_w)\cr
&=& -q\gots_w + \suml_{\ell(zs_i)>\ell(z)=\ell(w)} a_{w,z} \gots_{z}+
(1-q)\suml_{\ell(zs_i)>\ell(z)=\ell(w)}d_{w,z} \gots_z \cr
&=& -q\gots_w+ \suml_{\ell(zs_i)>\ell(z)=\ell(w)}[(1-q)d_{w,z}+a_{w,z}]\gots_z,
\end{eqnarray*}
where 
$a_{w,z}\in\{0, \pm 1, \pm q\}$, $d_{w,z}\in \bbz$, and
the sum extends over all $z\in S_n$ with $\ell(zs_i)>\ell(z)=\ell(w)$,

Substituting 
$$
b_{w,z}=\cases
{d_{w,z}, &\hbox{ if $a_{w,z}\not=\pm q$}\cr
d_{w,z}-a_{w,z}/q, &\hbox{ if $a_{w,z}=\pm q$}\cr}
$$
and
$$
c_{w,z}=\cases
{a_{w,z}, &\hbox{ if $a_{w,z}\not=\pm q$}\cr
a_{w,z}/q, &\hbox{ if $a_{w,z}=\pm q$}\cr}
$$
completes the proof.

\qed

\medskip

Imitating the proof of Theorem 4.1 we obtain

\medskip

\begin{note}
\noindent{\bf Theorem 6.5.}
$$
\chi_2^k(T_\mu)=\suml_{\ell(w)=k} \hbox{weight}_q^\mu(w) 
$$
where $\hbox{weight}_q^\mu(w)$ are the same weights as in Theorem 4.1.
\end{note}

\medskip

Combining Theorems 6.5 and 4.1 together with
Ram's result (Theorem 2.1) shows that

\smallskip

\begin{note}
\noindent{\bf Theorem 6.6.} The representation of ${\cal H}_n(q)$
induced by the $q$-commutators $A_i$ on the homogeneous components of 
the coinvariant algebra,
and the representation induced by
the transposed randomized operators $R_i$ on these components,
are equivalent.
\end{note}

\medskip

\begin{note}
\noindent{\bf Problem 6.7.} Calculate the coefficients $b_{w,z}$ in 
Proposition~6.4.
\end{note}

We conjecture that 
$b_{w,z}\in \{-1,0,1\}$. 

\bigskip

\section{\bf Final Remarks}

\medskip

\subsection{\bf Related Families of Operators} 

Consider the following family of $q$-commutators: 
$$
B_i := -[\de_i,X_{i+1}]_q \qquad(1\le i<n).
$$
This family is closely related to the $q$-commutators $A_i$.

\smallskip

\begin{note}
\noindent{\bf Fact 7.1.} The operators  $B_i$
satisfy the Hecke algebra relations (2.1)--(2.3).
\end{note}

\smallskip

\begin{note}
\noindent{\bf Proposition 7.2.} 
Let $D_i$ be operators on the polynomial ring $P_n$ of the form
$c_i + P_i(X_i,X_{i+1}) \de_i$, $1\le i<n$,
where $P_i$ are polynomials of two variables, and $c_i$ are constants.
If 
\begin{itemize}
\item[(1)] $D_i$ satisfy the Hecke algebra relations (2.1)--(2.3)
(with $q\not=0,-1$);
 
\item[(2)] $D_i$ are degree preserving (as operators on $P_n$);

\item[(3)] $\Lambda_n$, the subring of symmetric functions, is pointwise
invariant under all $D_i$, $1\le i<n$;
\end{itemize}
 then, for $n>2$, either $D_i =A_i$ ($\forall i$) or 
 $D_i =B_i$ ($\forall i$).
\end{note}

\medskip

Proposition 7.2 is related to a general theorem of
Lascoux and Sch\"utzen\-berger.

\smallskip

\begin{note}
\noindent{\bf LS Theorem.} [LS Theorem 1] 
Let $x_1,x_2,x_3$ be variables,  and let $s_i$ $i=1,2$ be the simple
transpositions as above. Let $D_i$, $i=1,2$ be linear operators on the ring of 
{\bf rational functions} $\bbc(x_1,x_2,x_3)$
(considered as a vector space over $\bbc$) defined by:
$$
D_i=P_i  +Q_i s_i
$$
where $P_i,Q_i\in \bbc(x_i,x_{i+1})$ are rational functions of
the corresponding pair of variables.

Assume that:
\begin{itemize}
\item[(1)] $D_1 D_2 D_1=D_2 D_1 D_2$;
\item[(2)] $D_1$ is invertible and $P_1\not= 0$.
\end{itemize}
Then:

$D_1$ and $D_2$ preserve the ring of {\bf polynomials}
$\bbc [x_1,x_2,x_3]$ if and only if there exist
$\alpha,\beta, \gamma,\delta,\eta\in \bbc$, so that
$\Delta:=\alpha\delta-\beta\gamma\not=0$, $\eta\not=0$,
$\eta\not=\Delta$, and
$$
P_i(x_i,x_{i+1})= (x_i-x_{i+1})^{-1}(\alpha x_i +\beta)
(\gamma x_{i+1}+\delta)
\hbox{ }\hbox{ }\hbox{ }
\hbox{ and } 
\hbox{ }\hbox{ }\hbox{ }
Q_i=\eta-P_i.
$$
Also, in that case, both $D_1$ and $D_2$ satisfy 
$D_i^2=\Delta D_i+\eta(\eta-\Delta)$.
\end{note}

\medskip 

Obviously, the initial conditions in this theorem are quite different
from those of Proposition 7.2; but the two families $A_i$ and $B_i$
are common solutions of both problems 
(for the LS theorem in the special case $\Delta=1-q$, $\eta=1$).
Note that for $\Delta=1-q$, $\eta=-q$ one gets two other families
of the $q$-commutator type, for which $\Lambda_n$ is $not$ pointwise
invariant.

It should be mentioned that the family
 $R_i$  of Sections 5 and 6 is $not$ obtainable from the $LS$ 
theorem (or from Proposition 7.2).

\medskip

\subsection{\bf Connections with Kazhdan-Lusztig Theory}

\smallskip

Theorem 3.3 has a remarkable analogue in Kazhdan-Lusztig theory.
In their seminal paper [KL] Kazhdan and Lusztig constructed a canonical basis
to Hecke algebra representations.
A basic  theorem in this theory describes the action of the generators $T_s$ 
on the canonical basis elements $C_w$. 

\smallskip

\begin{note}
\noindent{\bf Theorem 7.3.} [KL (2.3.a)--(2.3.c)] Let $W$ be a Coxeter group, $s$
  a Coxeter generator of $W$, 
$w\in W$, and $C_w$  the corresponding Kazhdan-Lusztig basis element. Then 
$$
T_s(C_w)=\cases{
-C_w,  &\hbox { if $\ell(sw)<\ell(w)$;}\cr
qC_w + q^\half \suml_{\ell(sz)<\ell(z)=\ell(w)} a_{w,z} C_z,
&\hbox { if $\ell(sw)>\ell(w)$,}\cr
}
$$
where the coefficients $a_{w,z}\in \bbz$ are independent of $q$.
\end{note}

\smallskip

This analogy leads to similar character formulas in the two theories;
see Theorems 2.2 and 4.1.
This surprising phenomenon seems to warrant further study.

\medskip

\subsection{\bf Probabilistic Aspects}

\smallskip

The parameter $q$ in the definition of the Hecke algebra may be interpreted 
as a {\it transition probability}. This gives a natural interpretation to 
the appearance of the coefficients $q$ and $1-q$ in the basic Hecke relation 
(2.3). 
This observation was fundamental to the definition of the randomized action 
in Section 5. The operators defined there interpolate between two well-studied
extreme cases: sorting ($q=0$) and mixing ($q=1$) by means of adjacent 
transpositions. They also form an interesting link between algebra and 
physics-motivated optimization.

\medskip

\subsection {\bf Other Weyl Groups}

Extension of all the above to other Weyl and Coxeter groups is highly desirable. Preliminary computations indicate that this may not be straightforward.

\end{document}